\documentclass[11pt]{amsart}          
\usepackage{amsmath,amsthm}
\usepackage{amssymb}            
\usepackage{verbatim}

\usepackage{array}
\usepackage{graphicx}
\usepackage [autostyle, english = american]{csquotes}
\MakeOuterQuote{"}

\usepackage{graphics}

\newcommand{\ml}{Martin-L\"{o}f }

\newcommand{\nat}{\in\omega}

\newcommand{\restr}{\!\!\restriction\!\!}
\newcommand{\st}{\;|\;}
\newcommand{\conv}{\!\!\downarrow}
\newcommand{\dvrg}{\!\!\uparrow}

\newcommand{\strA}{\mathcal{A}}
\newcommand{\strB}{\mathcal{B}}

\newcommand{\dega}{\mathbf{a}}
\newcommand{\degb}{\mathbf{b}}
\newcommand{\degd}{\mathbf{d}}
\newcommand{\degc}{\mathbf{c}}
\newcommand{\dego}{\mathbf{0}}

\newcommand{\ce}{c.e.\ }

\newtheorem{thm}{Theorem}[section]

\newtheorem{cor}[thm]{Corollary}
\newtheorem{lemma}[thm]{Lemma}

\theoremstyle{definition}  
\newtheorem{definition}[thm]{Definition}

\newtheorem{notation}[thm]{Notation}

\newtheorem*{construct}{Construction}

\newtheorem*{verify}{Verification}

\begin{document}

\title{Lowness for isomorphism, countable ideals, and computable traceability}

\author[Franklin]{Johanna N.Y.\ Franklin}
\address{Department of Mathematics \\ Room 306, Roosevelt Hall \\ Hofstra University \\ Hempstead, NY 11549-0114 \\ USA}
\email{johanna.n.franklin@hofstra.edu}
\thanks{The first author was supported in part by Simons Foundation Collaboration Grant \#420806. The authors thank Steffen Lempp for helpful comments during an early presentation 
of these results.}

\author[Solomon]{Reed Solomon}
\address{Department of Mathematics \\ University of Connecticut \\ Storrs, CT \\ USA}
\email{david.solomon@uconn.edu}

\date{\today}

\begin{abstract}
We show that every countable ideal of degrees that are low for isomorphism is contained in a principal ideal of degrees that are low for isomorphism by adapting an exact pair construction. We further show that within the hyperimmune-free degrees, lowness for isomorphism is entirely independent of computable traceability.
\end{abstract}

\maketitle

\section{Introduction}

Franklin and Solomon defined the Turing degrees that are low for isomorphism to be those degrees that are not useful for computing isomorphisms between computably presented structures: a degree $\mathbf{d}$ is low for isomorphism if, whenever it computes an isomorphism between computably presented structures $\strA$ and $\strB$, there is already a computable isomorphism between them \cite{fs-lowim}. This notion is robust in that the degrees that are low for isomorphism are precisely those that are low in various other contexts. Since graphs are universal structures \cite{hkss02}, it is clear that lowness for isomorphism for graphs corresponds to lowness for isomorphism.  Franklin and Turetsky have shown, more surprisingly, that it is also equivalent to lowness for paths: a degree is low for paths if, whenever it can compute an element of a  $\Pi^0_1$ class, either in Cantor space or Baire space, there must be a computable element in that class  \cite{ft-lowpaths}.

However, this class of Turing degrees has, so far, defied a full characterization that does not involve lowness for structures in some way. It differs from lowness classes like the low or superlow degrees or the $K$-trivials because no degree comparable to $\dego'$ (other than $\dego$) is low for isomorphism. 
It is known that the low for isomorphism degrees are not an ideal, that the 2-generic degrees form a subclass of them, and that no \ml random degree is low for isomorphism \cite{fs-lowim}. We can use \ml randomness as a boundary for the level of randomness that is incompatible with lowness for isomorphism \cite{fs-lowim}, and Franklin and Turetsky have constructed a 1-generic degree that is low for isomorphism and is not bounded by a 2-generic degree \cite{ft-1genlow}. In Section \ref{sec:HIF}, we will use the fact that any degree that can compute a separating set for two computably inseparable \ce sets is not low for isomorphism \cite{fs-lowim} and thus that no PA degree is low for isomorphism. However, none of this gives us the sought-after description. In fact, this class of Turing degrees is not a subclass of any standard degree-theoretic class of Turing degrees we have been able to identify: the hyperimmune degrees, the hyperimmune-free degrees, the minimal degrees, etc.

In this paper, we investigate the upward closure properties of the degrees that are low for isomorphism. As previously mentioned, these degrees are not closed under join: every 2-generic is low for isomorphism, and we can find two 2-generics whose join is above $\dego'$. We prove in Section \ref{sec:ExactPairs} that every increasing sequence of Turing degrees that is low for isomorphism has an exact pair whose join is low for isomorphism. It follows that every countable ideal within the low for isomorphism degrees is contained in a principal ideal within these degrees. Then, in Section \ref{sec:HIF}, we show that computable traceability is entirely independent of lowness for isomorphism inside the hyperimmune-free degrees. That is, within the hyperimmune-free degrees, there are degrees satisfying each of the four combinations of being computably traceable (or not) and being low
for isomorphism (or not). 

Our notation is largely standard. We use $\lambda$ to denote the empty string and $\sigma*\tau$ to represent the concatenation of the two binary strings $\sigma$ and $\tau$. 
For a tree $T \subseteq 2^{< \omega}$, we let $[T] = \{ f \mid \forall n \, ( f \restr n \in T )\}$ denote the set of infinite paths through $T$. For a pair of computable 
structures $\mathcal{A}$ and $\mathcal{B}$, we write $\mathcal{A} \cong_{\Delta^0_1} \mathcal{B}$ if there is a computable isomorphism $f: \mathcal{A} \rightarrow \mathcal{B}$. 

\section{Principal ideals}\label{sec:ExactPairs}

Recall that an exact pair $\degb$ and $\degc$ for an increasing sequence $\dega_0<_T\dega_1<_T\ldots$ is a pair of upper bounds for the sequence such that any degree below both $\degb$ and $\degc$ is below $\dega_n$ for some $n$. We begin by observing there are increasing sequences of degrees that are low for isomorphism and this result is therefore nonvacuous: in \cite{fs-lowim}, such a sequence was constructed using Mathias forcing.

\begin{thm}
\label{thm:exact}
Let $\dega_0<_T\dega_1<_T\ldots$ be an increasing sequence of Turing degrees that are low for isomorphism. This sequence has an exact pair $\degb$ and $\degc$ such that 
$\degb \oplus \degc$ is low for isomorphism, and hence $\degb$ and $\degc$ are low for isomorphism as well.
\end{thm}

The proof of this theorem is heavily based on the proof of the Kleene-Post-Spector result that every ascending sequence has an exact pair in \cite{oldsoare} (Chapter VI, Theorem 4.2). 

\begin{notation}
Let $A$ and $B$ be subsets of $\omega$, and let $f$ and $g$ be partial functions from $\omega$ to $2$.
We say that the \emph{$y$-section} of $A$ is 
$$A^{[y]}= \{\langle x,y\rangle \st \langle x,y\rangle\in A\},$$
and we define
$$A^{[<y]} = \bigcup \{ A^{[z]} \st z<y\}.$$
We further say that $A=^*B$ if the symmetric difference $(A-B)\cup(B-A)$ is finite.

Finally, we say that $f$ and $g$ are compatible and write $compat(f,g)$ if there is no $n$ such that $f(n)$ and $g(n)$ are defined and not equal.
\end{notation}

\begin{proof}
We begin by fixing a set $A_y\in \dega_y$ for each $y$ and defining $A = \{\langle x,y\rangle \st x\in A_y\}$. We fix an enumeration of the computably presented graphs $\langle \mathcal{G}_i\rangle$ with domain $\omega$. We will construct characteristic functions $f$ and $g$ of $B$ and $C$, respectively, in such a way as to guarantee that $B$ and $C$ form an exact pair for our sequence and that $B\oplus C$ is low for isomorphism. Our requirements are as follows:

\begin{description}
\item[$\mathcal{T}^B_y$] $B^{[y]} =^* A^{[y]}$.
\item[$\mathcal{T}^C_y$] $C^{[y]} =^* A^{[y]}$.
\item[$\mathcal{R}_{\langle e,i\rangle}$] If $\Phi^B_e = \Phi^C_i = h$ for some total $h$, then there is a $y$ such that $h\leq_T A^{[y]}$.
\item[$\mathcal{S}^{B\oplus C}_{\langle e,j,k\rangle}$] If $\Phi_e^{B\oplus C}$ is an isomorphism from $\mathcal{G}_j$ to $\mathcal{G}_k$, then there is a $y$ such that an isomorphism from $\mathcal{G}_j$ to $\mathcal{G}_k$ is computable in $A^{[y]}$.
\end{description}

The $\mathcal{T}$ requirements guarantee that $A_y$ is computable from $B^{[y]}$ and $C^{[y]}$ and thus from $B$ and $C$ for every $y$, so $B$ and $C$ will each be an upper bound for the $\dega_y$-sequence.  The $\mathcal{R}$ requirements are the exact pair requirements: anything that $B$ and $C$ can both compute must be computable from an $A_y$. These two types of requirements will be met as in the standard proof of the exact pair theorem. 
The $\mathcal{S}$ requirement ensures that $B \oplus C$ is low for isomorphism: if $B \oplus C$ can compute an isomorphism between two graphs, then one of the $A_y$ sets 
must be able to as well. Since each $A_y$ is low for isomorphism, there must therefore be a computable isomorphism between them.

\begin{construct}
At stage $s=0$, we set $f_0=g_0=\emptyset$.

At the beginning of stage $s+1$, we assume that $f_s$ and $g_s$ are fully defined on $\omega^{[<s]}$ and that the $\mathcal{T}$-requirements have been satisfied for all $y<s$; that is, 
$$(\forall y<s)[B_s^{[y]} =^* C_s^{[y]} =^* A^{[y]}].$$
We further assume that only finitely many values have been defined for $f_s$ and $g_s$ outside of $\omega^{[<s]}$. Below, we use $\sigma$ and $\tau$ as variables ranging over 
finite binary strings.

First we satisfy $\mathcal{R}_s$. Suppose that $s=\langle e,i\rangle$. We ask whether
$$(\exists \sigma, \tau)(\exists x)(\exists t)[(compat(\sigma,f_s)\wedge compat(\tau,g_s)) \wedge \Phi^\sigma_{e,t}(x)\conv \neq \Phi^\tau_{i,t}(x)\conv].$$
If such $\sigma$ and $\tau$ exist, we extend $f_s$ to $\hat{f}_s =f_s\cup \sigma$ and $g_s$ to $\hat{g}_s =g_s\cup\tau$; otherwise, we let $\hat{f}_s =f_s$ and $\hat{g}_s =g_s$.

Second, we satisfy $\mathcal{S}^{B \oplus C}_s$, where $s=\langle e,j,k\rangle$. We begin by trying to force $\Phi_e^{B \oplus C}$ to be partial. We ask whether there exist 
$x$, $\sigma$ and $\tau$ such that $\sigma$ is compatible with $\hat{f}_s$, $\tau$ is compatible with $\hat{g}_s$, and for every pair of strings $\hat{\sigma}$ and $\hat{\tau}$ 
for which $\hat{\sigma}$ extends $\sigma$ and is compatible with $\hat{f}_s$, and $\hat{\tau}$ extends $\tau$ and is compatible with $\hat{g}_s$, we have 
$\Phi_e^{\hat{\sigma} \oplus \hat{\tau}}(x) \dvrg$. If so, we extend $\hat{f}_s = \hat{f}_s \cup \sigma$ and $\hat{g}_s = \hat{g}_s \cup \tau$ and we are done with 
$\mathcal{S}^{B \oplus C}_{\langle e,j,k\rangle}$. 

If the answer is no, we try to force $\Phi_e^{B \oplus C}$ to fail to be onto $\mathcal{G}_k$. We ask whether there exist $y$, $\sigma$ and $\tau$ 
such that $\sigma$ is compatible with $\hat{f}_s$, $\tau$ is compatible with $\hat{g}_s$, and for every pair of strings $\hat{\sigma}$ and $\hat{\tau}$ for which $\hat{\sigma}$ 
extends $\sigma$ and is compatible with $\hat{f}_s$, and $\hat{\tau}$ extends $\tau$ and is compatible with $\hat{g}_s$, and for every $x$ such that 
$\Phi_e^{\hat{\sigma} \oplus \hat{\tau}}(x)\conv$, we have $\Phi_e^{\hat{\sigma} \oplus \hat{\tau}}(x) \neq y$. If so, we extend  
$\hat{f}_s = \hat{f}_s \cup \sigma$ and $\hat{g}_s = \hat{g}_s \cup \tau$ and we are again done with $\mathcal{S}^{B \oplus C}_{\langle e,j,k\rangle}$. 

If the answer is no again, we try to force $\Phi_e^{B \oplus C}$ to fail to be a partial isomorphism from some initial segment of $\mathcal{G}_j$ 
to $\mathcal{G}_k$. We ask whether there exist $x$, $\sigma$ and $\tau$ such that $\sigma$ is compatible with $\hat{f}_s$, $\tau$ is compatible with $\hat{g}_s$, 
$\Phi_e^{\sigma \oplus \tau}(u) \conv$ for all $u < x$, but $\Phi_e^{\sigma \oplus \tau} \restr x$ is not a partial isomorphism from $\mathcal{G}_j$ to 
$\mathcal{G}_k$. If the answer is yes, we extend $\hat{f}_s = \hat{f}_s \cup \sigma$ and $\hat{g}_s = \hat{g}_s \cup \tau$. If the answer is no, we leave $\hat{f}_s$ and 
$\hat{g}_s$ unchanged. In either case, we are done with $\mathcal{S}^{B \oplus C}_{\langle e,j,k\rangle}$.

Finally, we satisfy $\mathcal{T}^B_s$ by defining $f_{s+1}(x) = \hat{f}_s(x)$ for all $x \in \text{domain}(\hat{f}_s)$ and 
$f_{s+1}(x)=A(x)$ for all $x\in \omega^{[s]}-\text{domain}(\hat{f})$. Similarly, we satisfy $\mathcal{T}^C_s$ by defining $g_{s+1}(x) = \hat{g}_s(x)$ for all 
$x \in \text{domain}(\hat{g}_s)$ and $g_{s+1}(x)=A(x)$ for all $x\in \omega^{[s]}-\text{domain}(\hat{g})$. Because $\hat{f}_s$ and $\hat{g}_s$ are finite extensions of 
$f_s$ and $g_s$, we have maintained the inductive hypothesis that $f_{s+1}$ and $g_{s+1}$ are defined on all of $\omega^{[< s+1]}$ and are only defined on finitely many points 
outside of $\omega^{[< s+1]}$, and we have ensured that $$(\forall y<s+1)[B_{s+1}^{[y]} =^* C_{s+1}^{[y]} =^* A^{[y]}].$$
\end{construct}

\begin{verify}
We claim that $B$ and $C$ form an exact pair for the $\dega_y$-sequence and that $B \oplus C$ is low for isomorphism. It is immediate from the construction that 
$B$ and $C$ are upper bounds for this sequence because $A_s\equiv_T A^{[s]}\leq_T B^{[s]}\leq_T B$ and $A_s\equiv_T A^{[s]}\leq_T C^{[s]}\leq_T C$.

Now we persuade ourselves that if $B$ and $C$ both compute the same total function $h$, then $h$ must be computable from $A_y$ for some $y$. Fix the indices $e$ and $i$ such 
that $\Phi^B_e = \Phi^C_i=h$ and suppose that $\langle e,i\rangle=s$. If the answer to the $\mathcal{R}_s$ question were yes, then by the definition of $\hat{f}_s$ and $\hat{g}_s$, we would 
have $\Phi_e^B(x) = \Phi_e^{\sigma}(x) \neq \Phi_i^{\tau} = \Phi_i^C(x)$, contrary to hypothesis. Therefore, the answer to the $\mathcal{R}_s$ question is no. $\Phi_e^B$ is 
total, so for each $x$, there is a $\sigma$ compatible with $f_s$ such that $\Phi_e^B(x) \conv$, and similarly for $g_s$. It follows that for each $x$ 
and each $\sigma$ compatible with $f_s$ for which $\Phi^{\sigma}(x) \conv$, we have $\Phi^{\sigma}(x) = h(x)$. Therefore, the value of $h(x)$ can be computed from 
$f_s$ (and hence from $A_s$) by searching for any string $\sigma$ compatible with $f_s$ for which $\Phi_e^{\sigma}(x) \conv$. 

Finally, we argue that if $B \oplus C$ can compute an isomorphism between $\mathcal{G}_j$ and $\mathcal{G}_k$, then there is a computable isomorphism between the two 
graphs. Fix indices $e$, $j$ and $k$ such that $\Phi_e^{B \oplus C}$ is an isomorphism from $\mathcal{G}_j$ onto $\mathcal{G}_k$. Let $s=\langle e,j,k\rangle$. Because 
$\hat{f}_s \equiv_T \hat{g}_s \leq_T A_s$ and $A_s$ is low for isomorphism, it suffices to show that there is an isomorphism $h \leq_T \hat{f}_s \oplus \hat{g}_s$.  Consider the 
action for $\mathcal{S}^{B \oplus C}_{\langle e,j,k\rangle}$ at stage $s+1$. Since $\Phi_e^{B \oplus C}$ is a total isomorphism from $\mathcal{G}_j$ onto $\mathcal{G}_k$, the 
answer to the three questions for $\mathcal{S}^{B \oplus C}_{\langle e,j,k\rangle}$ must be no (else we would have taken a finite extension forcing this function not to be an 
isomorphism). We compute an isomorphism $h \leq_T \hat{f}_s \oplus \hat{g}_s$ by a back-and-forth argument as follows. 

Effectively in $\hat{f}_s \oplus \hat{g}_s$, we define sequences of finite strings $\sigma_t$ and $\tau_t$ which are compatible with $\hat{f}_s$ and $\hat{g}_s$, respectively, and a 
sequence of numbers $m_t$ such that $m_0 = 0$, $\sigma_0 = \tau_0 = \lambda$, $m_t < m_{t+1}$, $\sigma_{t+1}$ extends $\sigma_t$, $\tau_{t+1}$ extends 
$\tau_t$, and $\Phi_e^{\sigma_t \oplus \tau_t} \restr m_t$ is a partial isomorphism from $\mathcal{G}_j$ into 
$\mathcal{G}_k$ that (for $t >0$) includes $t-1$ in its range. Setting $h(t) = \Phi_e^{\sigma_{t+1} \oplus \tau_{t+1}}(t)$ gives the desired isomorphism. 

Assume $\sigma_t$, $\tau_t$ and $m_t$ have been defined. Let $h_t = \Phi_e^{\sigma_t \oplus \tau_t} \restr m_t$ be the partial isomorphism determined so far. 
Let $y$ be the least number not currently in the range of $h_t$. We define $\sigma_{t+1}$ and $\tau_{t+1}$ in two steps. First, we extend the range of our isomorphism to 
include $y$. Because the answer to the second question about $\mathcal{S}^{B \oplus C}_{\langle e,j,k\rangle}$ is no and because $\sigma_t$ and $\tau_t$ are compatible 
with $\hat{f}_s$ and $\hat{g}_s$ respectively, we can find (by searching effectively in $\hat{f}_s \oplus \hat{g}_s$) an $x$ and strings $\sigma_t'$ and $\tau_t'$ such that 
$\sigma_t'$ extends $\sigma_t$ and is compatible 
with $\hat{f}_s$, $\tau_t'$ extends $\tau_t$ and is compatible with $\hat{g}_s$ and $\Phi_e^{\sigma_t' \oplus \tau_t'}(x) = y$. Note that 
$x \geq m_t$ because $y$ was not in the domain of $h_t$. We now set $m_{t+1} = x+1$. 

Second, we make sure $h_{t+1}$ is defined on all the elements in $[m_t, m_{t+1})$. Since the answer to the first question about $\mathcal{S}^{B \oplus C}_{\langle e,j,k\rangle}$ 
is no and because $\sigma_t'$ and $\tau_t'$ are compatible with $\hat{f}_s$ and $\hat{g}_s$, respectively, we can find successive finite extensions of $\sigma_t'$ and $\tau_t'$ 
that are compatible with $\hat{f}_s$ and $\hat{g}_s$ and that force convergence on each value in $[m_t, m_{t+1})$ in turn. Let $\sigma_{t+1}$ and $\tau_{t+1}$ 
be the final such extensions. By construction, $\Phi_e^{\sigma_{t+1} \oplus \tau_{t+1}}(u) \conv$ for all $u < m_{t+1}$, and the range of 
$h_{t+1} = \Phi_e^{\sigma_{t+1} \oplus \tau_{t+1}} \restr m_{t+1}$ contains all $v \leq y$ (which, inductively, implies that the range contains $t$). Since the answer to the third 
question about $\mathcal{S}^{B \oplus C}_{\langle e,j,k\rangle}$ is no and since $\sigma_{t+1}$ and $\tau_{t+1}$ are compatible with $\hat{f}_s$ and $\hat{g}_s$, 
$h_{t+1}$ is a partial isomorphism. It is clear that 
the sequence of strings $\sigma_t$ and $\tau_t$ and numbers $m_t$ are uniformly computable in $\hat{f}_s \oplus \hat{g}_s$ and hence that the resulting isomorphism $h$ is 
computable from $\hat{f}_s \oplus \hat{g}_s$ as required. 
\end{verify}
\end{proof}

\begin{cor}
Every countable ideal within the low for isomorphism degrees is contained in a principal ideal within the low for isomorphism degrees.
\end{cor}

\begin{proof}
Let $\mathcal{I} = \{ \degd_i \mid i \in \omega \}$ be a countable ideal within the low for isomorphism degrees. For $j \in \omega$, define $\mathbf{e}_j = \oplus_{i \leq j} \degd_i$. 
By definition, $\mathbf{e}_0 \leq_T \mathbf{e}_1 \leq_T \ldots$ and $\dega \in \mathcal{I}$ if and only if $\dega \leq_T \mathbf{e}_j$ for some $j$. If this sequence is eventually 
constant, say equal to $\mathbf{e}$, then $\mathcal{I} = \{ \degd \mid \degd \leq_T \mathbf{e} \}$ and we are done. Otherwise, we can thin this sequence to a strictly increasing 
sequence $\dega_0 <_T \dega_1 <_T \ldots$ to which we can apply Theorem \ref{thm:exact}. In this case, $\mathcal{I} \subseteq \{ \degd \mid \degd \leq_T \degb \oplus \degc \}$.
\end{proof}

A similar construction can be used in the context of other downward closed sets of degrees such as the non-DNC degrees 
and the non-PA degrees. We recall their definitions here; for a summary of the properties of DNC and PA degrees, please see Sections 2.22 and 2.21 of \cite{dhbook}, respectively.

\begin{definition}
A Turing degree $\degd$ is DNC if it computes a total function $f$ such that for every $e$, $f(e)\neq \Phi_e(e)$, and it is PA if it is the degree of a complete extension of Peano 
arithmetic.
\end{definition}

The proof of Theorem \ref{thm:exact} does not require much modification in either the non-DNC or non-PA case. Suppose that we have an increasing sequence 
$\dega_0<_T\dega_1<_T\ldots$ of degrees that are not DNC (or not PA). The $\mathcal{T}^B_y$, $\mathcal{T}^C_y$, and $\mathcal{R}_{\langle e,i\rangle}$ requirements are 
stated and satisfied as before.  For the non-DNC case, the $\mathcal{S}^{B \oplus C}$ requirements become
\begin{description}
\item[$\mathcal{S}^{B \oplus C}_e$] If $\Phi^{B \oplus C}_e$ is total, then there is an $n$ such that $\Phi^{B \oplus C}_e(n)=\varphi_n(n)$.
\end{description}
To satisfy $\mathcal{S}^{B \oplus C}_e$, assume $s=e+1$ and we have handled the appropriate $\mathcal{R}$ requirement to obtain $\hat{f}_s$. 
We first try to force $\Phi_e^{B \oplus C}$ 
to be partial by asking the same question as in the proof of Theorem \ref{thm:exact}. If the answer is yes, we guarantee partiality by extending 
$\hat{f}_s = \hat{f}_s \cup \sigma$ 
and $\hat{g}_s = \hat{g}_s \cup \tau$, and we are done with $\mathcal{S}^{B \oplus C}_e$. If we cannot force partiality, we try to force $\Phi_e^{B \oplus C}$ to be 
non-DNC by asking if there exist $n$, $\sigma$ and $\tau$ such that $\varphi_n(n) \conv$, $\sigma$ is compatible with $\hat{f}_s$, $\tau$ is compatible 
with $\hat{g}_s$ and $\Phi_e^{\sigma \oplus \tau}(n) \conv = \varphi_n(n)$. Again, if the answer is yes, we extend $\hat{f}_s = \hat{f}_s \cup \sigma$ 
and $\hat{g}_s = \hat{g}_s \cup \tau$. If the answer is no, we leave $\hat{f}_s$ and $\hat{g}_s$ unchanged.

To verify that this works, note that if the answer to either question is yes, then we have forced $\Phi_e^{B \oplus C}$ to be partial or non-DNC and hence have satisfied 
$\mathcal{S}^{B \oplus C}_e$. Assume the answers are both no. In this case, we claim that there is a total function $h \leq_T \hat{f}_s \oplus \hat{g}_s$ such that for every $n$, 
$\varphi_n(n) \neq h(n)$. Since $A_s \geq_T \hat{f}_s \equiv_T \hat{g}_s$, this claim implies $A_s$ is DNC, contradicting our hypothesis. 

To prove the claim, we define $h$ as follows. Effectively in $\hat{f}_s \oplus \hat{g}_s$, we construct sequences of finite strings $\sigma_t$ and $\tau_t$ which are 
compatible with $\hat{f}_s$ and $\hat{g}_s$ respectively such that $\sigma_0 = \tau_0 = \lambda$, $\sigma_{t+1}$ 
extends $\sigma_t$, $\tau_{t+1}$ extends $\tau_t$, and for $u \leq t$, $\Phi_e^{\sigma_{t+1} \oplus \tau_{t+1}}(u) \conv \neq \varphi_u(u)$. Defining 
$h(t) = \Phi_e^{\sigma_{t+1} \oplus \tau_{t+1}}(t)$ gives the desired DNC function. 

To define the sequences of finite strings, assume $\sigma_t$ and $\tau_t$ have been defined with the desired properties. Because the answer to the first question for 
$\mathcal{S}^{B \oplus C}_e$ is no, we can find (by searching 
effectively in $\hat{f}_s$ and $\hat{g}_s$) strings $\sigma_{t+1}$ extending $\sigma_t$ and $\tau_{t+1}$ extending $\tau_t$ which are compatible with 
$\hat{f}_s$ and $\hat{g}_s$, respectively, such that $\Phi_e^{\sigma_{t+1} \oplus \tau_{t+1}}(t) \conv$. Since the answer to the second question about 
$\mathcal{S}^{B \oplus C}_e$ is no, we must have $\Phi_e^{\sigma_{t+1} \oplus \tau_{t+1}}(t) \neq \varphi_t(t)$ as required to complete the non-DNC construction.

The non-PA case is essentially the same, using the fact (from \cite{js72-2}) that a degree is PA if and only if it computes a $\{ 0,1 \}$-valued DNC function. 
The $\mathcal{S}^{B \oplus C}$ requirements become
\begin{description}
\item[$\mathcal{S}^{B \oplus C}_e$] If $\Phi^{B \oplus C}_e$ is total and $\{0,1\}$-valued, then there is an $n$ such that $\Phi^{B \oplus C}_e(n)=\varphi_n(n)$.
\end{description}
To meet $\mathcal{S}^{B \oplus C}_e$, we first try to force $\Phi^{B \oplus C}_e$ to be partial. If that fails, we try to force $\Phi^{B \oplus C}_e(n) > 1$ for some 
$n$. If that fails as well, we try to force $\Phi^{B \oplus C}_e(n) = \varphi_n(n)$ for some $n$. In the case when we cannot force any of these conditions, we use 
$\hat{f}_s \oplus \hat{g}_s$ to compute a $\{0,1\}$-valued DNC function as above. 

These arguments give us the following result. Once again, we note that this is not vacuous: every Turing incomplete \ce degree is both non-PA \cite{js72-2} and non-DNC \cite{arslanov81,jlss89}.

\begin{cor}
Every countable ideal within the non-DNC degrees is contained in a principal ideal within the non-DNC degrees, and every countable ideal within the non-PA degrees is contained in a principal ideal within the non-PA degrees.
\end{cor}

We end this section by noting that there are also downward closed classes of degrees which do not satisfy an exact pair theorem of the type considered here. For example, it is 
straightforward to construct an ideal within the low degrees which is not contained in any principal ideal topped by a low degree. Another related result (involving a 
class of degrees that is not downward closed) was proved by Ambos-Spies: there is a uniformly c.e.~ascending sequence of degrees $\dega_0<_T\dega_1<_T\ldots$ 
which has no exact pair of c.e.~degrees. (See Chapter IX, Exercise 4.4 in \cite{oldsoare}.)

\section{Hyperimmune-free degrees}\label{sec:HIF}

We now turn our attention to the hyperimmune-free degrees. A degree $\degd$ is hyperimmune-free if for all functions 
$f \leq_T \degd$, there is a computable function $g$ such that $f(n)\leq g(n)$ for every $n$. Terwijn and Zambella defined a natural subclass of the hyperimmune-free degrees, the computably traceable degrees. In the definition below, an order function is a nondecreasing and unbounded computable function from $\omega$ to $\omega$.

\begin{definition}[Terwijn and Zambella \cite{tz01}]
A set $A$ is computably traceable if there is an order function $p$ such that for every $f\leq_T A$, there is a computable function $r$ such that for every $n\nat$, 
\begin{enumerate}
\item $f(n)\in D_{r(n)}$ and
\item $|D_{r(n)}|\leq p(n)$.
\end{enumerate}  
A degree is computably traceable if it contains a computably traceable set. 
\end{definition}

Note that every computably traceable degree is hyperimmune free.
Terwijn and Zambella proved that if there is an order function $p$ as in this definition, then any order function will actually work and that the degrees that are 
low for Schnorr tests are precisely the computably traceable degrees \cite{tz01}. It follows that a \ml random degree cannot be computably traceable, because 
a \ml random degree is also Schnorr random, and hence by Terwijn and Zambella's characterization is not computably traceable.

It is tempting to think that the computably traceable degrees could be precisely the hyperimmune-free degrees that are low for isomorphism. $\dego$ is an obvious example 
of a computably traceable degree that is low for isomorphism, and since every \ml random degree is neither computably traceable nor low for isomorphism (by 
Theorem 4.1 in \cite{fs-lowim}), there are hyperimmune-free degrees that are neither computably traceable nor low for isomorphism. However, in this section, we prove that 
the hyperimmune-free degrees that are low for isomorphism are neither a subclass nor a superclass of the computably traceable degrees in Theorems 
\ref{thm:CompTrNotLowIM} and \ref{thm:NotCompTrLowIM}.

\begin{thm}\label{thm:CompTrNotLowIM}
There is a computably traceable degree that is not low for isomorphism.
\end{thm}

\begin{proof}
We build a computably traceable set $C$ which is not low for isomorphism. To make $C$ not low for isomorphism, we build a pair of 
computably inseparable c.e.~sets $A$ and $B$, and we ensure that $C$ is a separating set for this pair. We denote the class of separating sets by $\text{Sep}(A,B)$. 
This method suffices because any set which computes a separating set for any pair of computably inseparable c.e.~sets is not low for isomorphism by Theorem 3.4 in \cite{fs-lowim}. 

Formally, we build a $\{ 0,1 \}$-valued partial computable function 
$\psi$ in stages as $\psi_s$ with $\psi_s \subseteq \psi_{s+1}$ and $\psi = \cup_s \psi_s$.  We set $A = \{ n \mid \psi(n) = 1 \}$, $B = \{ n \mid \psi(n) = 0 \}$,  
$A_s = \{ n \mid \psi_s(n) = 1 \}$, and $B_s = \{ n \mid \psi_s(n) = 0 \}$ and define a sequence of trees such that $[V_s] = \text{Sep}(A_s,B_s)$ by
\[
V_s = \{ \sigma \in 2^{< \omega} \mid \forall n \leq |\sigma| \, \left( \psi_s(n) \! \downarrow \rightarrow \psi_s(n) = \sigma(n) \right) \}.
\]
These trees form a uniformly computable nested sequence $V_0 \supseteq V_1 \supseteq \ldots$ such that $V = \cap_s V_s$ is a $\Pi^0_1$ tree for which 
$[V] = \text{Sep}(A,B)$. To make the domain of $\psi$ coinfinite, we use movable markers $\delta_s(i)$ to denote the $i$-th element of the complement of $A_s \cup B_s$, 
and we ensure that $\delta(i) = \lim_s \delta_s(i)$ exists for each $i$. Note that by definition, $\psi_s(\delta_s(i)) \! \uparrow$. 

To make $A$ and $B$ computably inseparable, we need to meet the requirements
\[
\mathcal{R}_e: \Phi_e \text{ is not the characteristic function for a separating set of } A \text{ and } B. 
\]
We say that $\mathcal{R}_e$ is \textit{met at stage} $s+1$ if there is an $x$ such that $\Phi_{e,s}(x) \! \downarrow \neq \psi_s(x) \downarrow$. We say that $\mathcal{R}_e$ is \textit{eligible to act at 
stage} $s+1$ if $\mathcal{R}_e$ is not met at stage $s+1$ and $\Phi_{e,s}(\delta_s(e)) \! \downarrow$. In this case, we can act to meet $\mathcal{R}_e$ by defining 
$\psi_{s+1}(\delta_s(e)) = 1 - \Phi_{e,s}(\delta_s(e))$. Since $\psi_s(\delta_s(e)) \! \uparrow$, this definition does not conflict with a previously defined value of $\psi$. 

Our set $C$ will be an infinite path in $V$. To make $C$ computably traceable, we meet the requirements
\[
\mathcal{S}_e: \Phi_e^C \text{ total } \Rightarrow (\exists  \text{ computable } r_e)(\forall n)\left[ \Phi_e^C(n) \in D_{r_e(n)} \wedge |D_{r_e(n)}| \leq 2^n \right].
\]
To meet $\mathcal{S}_e$, we ensure that for all $X \in [V]$, if $\Phi_e^X$ is total, then for all $n$, 
$\Phi_e^{X \upharpoonright \delta(e+n)}(n) \! \downarrow$. To describe the general strategy, let $\tau_{k,0}^s, \ldots, \tau_{k,2^k-1}^s$ denote the nodes in $V_s$ of length 
$\delta_s(k)$ ordered lexicographically (these are the branching nodes on the $k$-th branching level of $V_s$), and let $\tau_{k,\ell} = \lim_s \tau_{k,\ell}^s$ denote the limit 
(which we will prove exists below). Fix $k$. For each $e \leq k$ and each $\ell < 2^k$, we attempt to define 
$\psi$ so that $\Phi_e^{\tau_{k,\ell}}(k-e) \! \downarrow$ if possible. After the construction of $V$, we choose a particular path $C \in [V]$ that makes 
$\Phi_e^C$ partial if possible, and if this is impossible, we prove that an appropriate computable function $r_e(n)$ exists. 

Before giving the formal construction, we make the following definitions. Let $s$ be a stage, $e$ be an index, $k \geq e$, and $\ell < 2^k$. 
$\mathcal{S}_e$ is \textit{primed to act at} $\tau_{k,\ell}^s$ if for all $j$ such that $e \leq j < k$, $\Phi_{e,s}^{\tau_{k,\ell}^s \upharpoonright \delta_s(j)}(j-e) \downarrow$. 
We say that $\mathcal{S}_e$ is 
\textit{eligible to act at} $\tau_{k,\ell}^s$ if it is primed to act at $\tau_{k,\ell}^s$ and there is a $\sigma \in V_s$ such that $\tau_{k,l}^s \subseteq \sigma$ and 
$\Phi_{e,s}^{\sigma}(k-e) \! \downarrow$. In this case, we say $\mathcal{S}_e$ is eligible to act at $\tau_{k,\ell}^s$ \textit{via} $\sigma$. 
We say that $\mathcal{S}_e$ is \textit{currently satisfied at} $\tau_{k,\ell}^s$ if either $\mathcal{S}_e$ is not primed to act at $\tau_{k,\ell}^s$ or 
$\Phi_{e,s}^{\tau_{k,\ell}^s}(k-e) \! \downarrow$. 

We now give the construction of $\psi$. At stage $0$, we set $\psi_0 = \emptyset$. At an even stage $s+1$, we check whether there is an index $e \leq s$ such that 
$\mathcal{R}_e$ is eligible to act. If not, set $\psi_{s+1} = \psi_s$ and end the stage. If so, let $e$ be the least such index and define $\psi_{s+1}$ by $\psi_{s+1}(n) = \psi_s(n)$ 
for all $n \in \text{domain}(\psi_s)$ and $\psi_{s+1}(\delta_s(e)) = 1- \Phi_e(\delta_s(e))$. It follows that $\delta_{s+1}(i) = \delta_s(i)$ for $i < e$ and 
$\delta_{s+1}(i) = \delta_s(i+1)$ for $i \geq e$. 

At an odd stage $s+1$, check whether there is a level $k < s$ for which there is a node $\tau_{k,\ell}^s$ with $\ell < 2^k$ and an index $e \leq k$ such that $\mathcal{S}_e$ is eligible 
to act at $\tau_{k,\ell}^s$ and $\mathcal{S}_e$ is not currently satisfied at $\tau_{k,\ell}^s$. If not, set $\psi_{s+1} = \psi_s$ and end the stage. If so, fix the least such level $k$, then the least 
corresponding index $e < k$, and finally the least corresponding $\ell < 2^k$.  Fix $\sigma \in V_s$ such that 
$\tau_{k,\ell}^s \subseteq \sigma$ and $\Phi_{e,s}^{\sigma}(k-e) \! \downarrow$. Because $\mathcal{S}_e$ is not currently satisfied at $\sigma$, we know that $\tau_{k,\ell}^s \neq \sigma$. 
Define $\psi_{s+1}$ as follows and then end the stage. For each $n \in \text{domain}(\psi_s)$, set $\psi_{s+1}(n) = \psi_s(n)$. For each $n$ such that 
$|\tau_{k,\ell}^s| \leq n < |\sigma|$, set $\psi_{s+1}(n) = \sigma(n)$. (Because $\sigma \in V_s$, there is no conflict between these two cases.) 
Since $\tau_{k,\ell}^s \neq \sigma$, we know there is an $m > 0$ such that 
$\delta_s(k+j) \in \text{domain}(\psi_{s+1})$ for all $j < m$. In particular, $\delta_{s+1}(i) = \delta_s(i)$ for $i < k$, and $\delta_{s+1}(i) = \delta_s(i+m)$ for 
$i \geq k$. Furthermore, it follows that $\sigma = \tau_{k,\ell}^{s+1} \in V_{s+1}$ and $|\sigma| = \delta_{s+1}(k)$. 

This completes the construction of the tree $V$ and the sets $A$ and $B$. After verifying the necessary properties of this construction, we will choose a path $C$ 
through $V$ such that the degree of $C$ is computably traceable but not low for isomorphism. It is straightforward to verify the following properties of the construction of $V$. 
\begin{enumerate}
\item[(V1)] Each $\mathcal{R}_e$ requirement acts at most once. 
\item[(V2)] The only requirements that can cause $\delta_{s+1}(k) \neq \delta_s(k)$ are $\mathcal{R}_e$ and $\mathcal{S}_e$ with $e \leq k$. 
\item[(V3)] Suppose that at stage $s+1$, $\mathcal{S}_e$ acts at $\tau_{k,\ell}^s$ with $\sigma$. 
\begin{itemize}
\item $\tau_{k,\ell}^{s+1} = \sigma$ and therefore $\Phi_{e,s}^{\tau_{k,\ell}^{s+1}}(k-e) \! \downarrow$. 
\item Writing $\sigma = \tau_{k,\ell}^s * \mu$, we have that for all $j < 2^k$, $\tau_{k,j}^{s+1} = \tau_{k,j}^s*\mu$. 
\item For all $i < k$ and all $j < 2^k$, if $\Phi_{i,s}^{\tau_{k,j}^s}(k-i) \! \downarrow$, then $\Phi_{i,s}^{\tau_{k,j}^{s+1}}(k-i) \! \downarrow$.
\end{itemize}
\end{enumerate}

\begin{lemma}
For each $k$, $\lim_s \delta_s(k) = \delta(k)$ exists, and for each $\ell < 2^k$, $\tau_{k,\ell} = \lim_s \tau_{k,\ell}^s$ exists.
\end{lemma}

\begin{proof}
We proceed by induction on $k$. For $k=0$, let $t$ be a stage such that $\mathcal{R}_0$ never acts after stage $t$. For $s \geq t$, $\tau_{0,0}^s$ remains constant unless $\mathcal{S}_0$ acts at 
$\tau_{0,0}^s$. By construction, if $\mathcal{S}_0$ is eligible to act at $\tau_{0,0}^s$ and is not currently satisfied at $\tau_{0,0}^s$, then $\mathcal{S}_0$ will act at $\tau_{0,0}^s$ by defining 
$\tau_{0,0}^{s+1}$ to make $\Phi_{0,s}^{\tau_{0,0}^{s+1}}(0) \! \downarrow$. $\mathcal{S}_0$ will be  currently satisfied at $\tau_{0,0}^{s+1}$ and will not act again.

For the inductive case of $k > 0$, fix a stage $t$ such that $\delta_t(j)$ has reached its limit for all $j < k$ and such that $\mathcal{R}_k$ does not act after stage $t$. It follows that the nodes 
$\tau_{j,\ell}^t = \tau_{j,\ell}$ for $j < k$ and $\ell < 2^j$ have reached their limiting values. In particular, no requirement $\mathcal{R}_j$ for $j \leq k$ acts after stage $t$, 
and no requirement $\mathcal{S}_e$ acts at any node $\tau_{j,\ell}$ after stage $t$. Therefore, the only requirements that could cause $\delta_{s+1}(k) \neq \delta_s(k)$ are $\mathcal{S}_e$ 
requirements with $e \leq k$ acting at $\tau_{k,\ell}^s$ for some $\ell < 2^k$. By construction, these requirements have the highest priority to act at any odd stage.

It suffices to show that for each $e \leq k$ and $\ell < 2^k$, the requirement $\mathcal{S}_e$ can act at $\tau_{k,\ell}^s$ during at most one stage $s+1 > t$. Suppose 
$\mathcal{S}_e$ acts at $\tau_{k,\ell}^s$ during stage $s+1$. By (V3), we have $\Phi_{e,s}^{\tau_{k,\ell}^{s+1}}(k-e) \! \downarrow$, and this convergence is preserved whenever 
another requirement $\mathcal{S}_i$ with $i \leq k$ acts at a node $\tau_{k,j}^u$ at a future stage. Therefore, at every future stage $u$, $\mathcal{S}_e$ will remain currently satisfied 
at $\tau_{k,\ell}^u$ and hence will never act again at this node. 
\end{proof}

\begin{lemma}
Each $\mathcal{R}_e$ requirement is satisfied, so $A$ and $B$ are recursively inseparable.
\end{lemma}

\begin{proof}
Fix $e$ and let $t$ be a stage at which the value of $\delta_t(e)$ has reached its limit. If $\mathcal{R}_e$ is met at stage $t$, then $\mathcal{R}_e$ has already been permanently satisfied. 
If not, then we claim $\Phi_e(\delta(e)) \! \uparrow$ and hence $\mathcal{R}_e$ is satisfied. Suppose $\Phi_{e,s}(\delta(e)) \! \downarrow$ at some odd stage $s+1>t$. By construction, $\mathcal{R}_e$ 
will act at stage $s$ unless some higher priority $\mathcal{R}_i$ is also eligible to act. In either case, we have $\delta_{s+1}(e) \neq \delta_s(e)$, contradicting the choice of $t$. 
\end{proof}

\begin{lemma}
\label{lem:converge}
Let $X \in [V]$ be such that $\Phi_e^X$ is total. Then for all $n$, $\Phi_e^{X \upharpoonright \delta(n+e)}(n) \! \downarrow$. 
\end{lemma}

\begin{proof}
For a contradiction, fix the least $n$ for which $\Phi_e^{X \upharpoonright \delta(n+e)}(n) \! \uparrow$. Let $k=n+e$, so $n=k-e$. Let $t$ be a stage such that the 
$\delta(i)$s have reached their limits for $i \leq k$ and $\Phi_{e,t}^{X \upharpoonright \delta(j)}(j-e) \! \downarrow$ for all $j$ such that $e \leq j < k$. 
Fix $\ell < 2^k$ such that $X \upharpoonright \delta(k) = \tau_{k,\ell}^t = \tau_{k, \ell}$. At all stages $u \geq t$, $\mathcal{S}_e$ is primed to act at $\tau^u_{k, \ell} = \tau_{k,\ell}$, 
and therefore, since $\Phi_e^{\tau_{k,\ell}}(k-e) \! \uparrow$, $\mathcal{S}_e$ is not currently satisfied at $\tau^u_{k,\ell}$.

Fix $\sigma \subseteq X$ such that $\Phi_e^{\sigma}(k-e) \! \downarrow$ and note that $\tau_{k,\ell} \subseteq \sigma \in V$.  For any odd stage $u+1>t$ such that 
$\Phi_{e,u}^{\sigma}(k-e) \! \downarrow$, $\mathcal{S}_e$ is eligible to act at $\tau^u_{k,\ell}$ with $\sigma$, and therefore, we would define $\tau_{k,\ell}^{u+1} = \sigma$ for a 
contradiction.
\end{proof}

We can now define the desired path $C \in [V]$. For $\tau_{k,\ell}$, let $\text{Succ}_0(\tau_{k,\ell}) = \tau_{k+1,\ell'}$ where $\tau_{k,\ell}*0 \subseteq \tau_{k+1,\ell'}$. That is, by 
definition, $\tau_{k,\ell}$ is a branching node in $V$, and we let $\text{Succ}_0(\tau_{k,\ell})$ be the first branching node extending 
$\tau_{k,\ell}*0$. For us, the important points will be that $\tau_{k,\ell} \subseteq \text{Succ}_0(\tau_{k,\ell})$ and $|\text{Succ}_0(\tau_{k,\ell})| = \delta(k+1)$, so 
$\text{Succ}_0(\tau_{k,\ell}) \restr \delta(k) = \tau_{k,\ell}$. 

We define $C$ by specifying its initial segments $\sigma_0 \subseteq \sigma_1 \subseteq \ldots$. For each $e$, we will have $\sigma_e = \tau_{k,\ell}$ for some $k \geq e+1$ and 
$\ell < 2^k$. To define $\sigma_0$, if there exist $k$ and $\ell < 2^k$ such that $\Phi_0^{\tau_{k,\ell}}(k) \! \uparrow$, then fix the least such pair and set 
$\sigma_0 = \text{Succ}_0(\tau_{k,\ell})$. Otherwise, set $\sigma_0 = \text{Succ}_0(\tau_{0,0}) =  \tau_{1,0}$. 

Assume $\sigma_{e-1} = \tau_{a,b}$ has been defined with $e \leq a$ and $b < 2^a$. We consider the three possible scenarios separately. First, if there is an $n < a-e$ such that 
$\Phi_e^{\sigma_{e-1}}(n) \uparrow$, then set $\sigma_e = \text{Succ}_0(\sigma_{e-1})$. Second, suppose there is no such $n$ but there is a $k \geq a$ and 
$\ell < 2^k$ such that $\sigma_{e-1} \subseteq \tau_{k,\ell}$ and $\Phi_e^{\tau_{k,\ell}}(k-e) \! \uparrow$. In this case, fix the least such pair and set 
$\sigma_e = \text{Succ}_0(\tau_{k,\ell})$. Third, if neither of these cases applies, set $\sigma_e = \text{Succ}_0(\sigma_{e-1})$. 

It is immediate that each $\sigma_e$ has the form $\tau_{k,\ell}$ for some $k \geq e+1$ and that $C \in [V]$. Therefore, the degree of $C$ is not low for isomorphism. 
It remains for us to verify that each $\mathcal{S}_e$ requirement is met.

\begin{lemma}
Each requirement $\mathcal{S}_e$ is met. 
\end{lemma} 

\begin{proof}
We proceed by induction on $e$. For $e=0$, consider how $\sigma_0$ was defined. First, suppose $\sigma_0 = \text{Succ}_0(\tau_{k,\ell})$ where 
$\Phi_0^{\tau_{k,\ell}}(k) \! \uparrow$. In this case, we claim that $\Phi_e^C$ is not total. If $\Phi_0^C$ is total, then by 
Lemma \ref{lem:converge}, then $\Phi_0^{C \upharpoonright \delta(k)}(k) \! \downarrow$. But $C \! \upharpoonright \! \delta(k) = \sigma_0 \! \upharpoonright \! \delta(k) = \tau_{k,\ell}$, 
giving a contradiction. 

Second, suppose that $\sigma_0 = \text{Succ}_0(\tau_{0,0})$ because for every $k$ and $\ell < 2^k$, $\Phi_0^{\tau_{k,\ell}}(k) \! \downarrow$. We define a 
computable function $r_0(k)$ such that for all $k$, $|D_{r_0(k)}| \leq 2^k$ and $\Phi_0^C(k) \in D_{r_0(k)}$. For each $k$, run the construction of $V_s$ until a stage  
$s$ such that $\Phi_{0,s}^{\tau_{k,\ell}^s}(k) \! \downarrow$ for all $\ell < 2^k$. We must eventually see this convergence because there is a stage after which 
$\tau^s_{k,\ell} = \tau_{k,\ell}$ for all $\ell < 2^k$, and we know $\Phi_0^{\tau_{k,\ell}}(k) \! \downarrow$.  
Once we find such a stage $s$, let $D_{r_0(k)} = \{ \Phi_{0,s}^{\tau^s_{k,\ell}}(k) \mid \ell < 2^k \}$. Clearly, we have 
$|D_{r_0(k)}| \leq 2^k$ as required. Furthermore, for all $X \in [V_s]$, we have $\Phi_0^X(k) \! \downarrow \in D_{r_0(k)}$. Since $C \in [V] \subseteq [V_s]$, we have 
$\Phi_0^C(k) \in D_{r_0(k)}$ as required. 

The verification for $e > 0$ is similar. Assume $\sigma_{e-1} = \tau_{a,b}$ and consider how $\sigma_e$ was defined. First, suppose there is an $n < a-e$ such that 
$\Phi_e^{\sigma_{e-1}}(n) \! \uparrow$. We claim that $\Phi_e^C$ is not total. Suppose for a contradiction that $\Phi_e^C$ is total. Fix $j < a$ such that 
$n = j-e$, so $n+e=j$. By Lemma \ref{lem:converge}, $\Phi_e^{C \upharpoonright \delta(j)}(n) \! \downarrow$. However, 
$\sigma_{e-1} = \tau_{a,b} \subseteq C$ and $j < a$ imply that $C \! \upharpoonright \! \delta(j) = \sigma_{e-1} \! \upharpoonright \! \delta(j)$ and, therefore, 
$\Phi_e^{\sigma_{e-1}}(n) \! \downarrow$ for a contradiction. 

Second, suppose there is no such $n$, but there is a $k \geq a$ and $\ell < 2^k$ such that $\sigma_{e-1} \subseteq \tau_{k,\ell}$ and $\Phi_e^{\tau_{k,\ell}}(k-e) \! \uparrow$. 
Fix the $k$ and $\ell$ such that $\sigma_e = \text{Succ}_0(\tau_{k,\ell})$. We claim that $\Phi_e^C$ is not total. If it is 
total, then by Lemma \ref{lem:converge}, $\Phi_e^{C \upharpoonright \delta(k)}(k-e) \! \downarrow$. However, $\tau_{k,\ell} \subseteq C$, so 
$C \! \upharpoonright \! \delta(k) = \tau_{k,\ell}$ and hence $\Phi_e^{\tau_{k,\ell}}(k-e) \! \downarrow$ for a contradiction. 

Finally, suppose that we are not in either of these cases. We define a computable function $r_e(k)$ such that for all $k$, $|D_{r_e(k)}| \leq 2^k$ and $\Phi_e^C(k) \in D_{r_e(k)}$. 
Fix a stage $t$ such that $\tau^t_{a,b}$ has reached its limit $\tau_{a,b} = \sigma_{e-1}$. For $n < a-e$, set $D_{r_e(n)} = \{ \Phi_e^C(n) \} = \{ \Phi_e^{\sigma_{e-1}}(n) \}$.

For $n \geq a-e$, let $k = n+e$, so $k-e=n$. Run the construction of $V_s$ for $s > t$ until we see 
an $s$ such that $\Phi_{e,s}^{\tau_{k,\ell}^s}(k-e) \! \downarrow$ for all $\ell < 2^k$ such that $\tau_{a,b} \subseteq \tau_{k,\ell}$. 
As before, we must eventually see this convergence since there is a stage $s$ such that $\tau^s_{k,\ell} = \tau_{k,\ell}$ for all $\ell < 2^k$ with $\tau_{a,b} \subseteq \tau_{k,\ell}$, 
and by the case assumption, $\Phi_{e}^{\tau_{k,\ell}}(k-e) \! \downarrow$.  Define 
\[
D_{r_e(n)} = D_{r_e(k-e)} = \{ \Phi_{e,s}^{\tau^s_{k,\ell}}(k-e) \mid \ell < 2^k \text{ with } \tau_{a,b} \subseteq \tau^s_{k,\ell} \}. 
\]

We claim that for all $n \geq a-e$, we have $\Phi_e^C(n) \in D_{r_e(n)}$. Fix $n \geq a-e$ and let $k = n+e$. Let $X \in [V_s]$ be arbitrary with $\tau_{a,b} \subseteq X$. 
Let $\ell < 2^k$ be such that $\tau^s_{k,\ell} \subseteq X$. Since $k = n+e \geq a$, we have $\tau_{a,b} \subseteq \tau^s_{k, \ell}$ and therefore, 
$\Phi_e^X(n) = \Phi_e^X(k-e) = \Phi^{\tau^s_{k,\ell}}(k-e) \downarrow \in D_{r_e(n)}$. Since $X \in [V_s]$ was arbitrary and since $C \in [V] \subseteq [V_s]$, we have 
$\Phi_e^C(n) \in D_{r_e(n)}$ as required. 

Finally, we claim that $|D_{r_e(n)}| \leq 2^n$. This is clear for $n < a-e$ since $|D_{r_e(n)}| = 1$ in this case. Assume $n \geq a-e$ and let $k$ be such that 
$k = n+e \geq a$. Note that 
\[
|D_{r_e(n)}| \leq |\{ \ell < 2^k : \tau_{a,b} \subseteq \tau_{k,\ell}^s \}| \leq 2^{k-a} = 2^{n+e-a}.
\]
Since $a \geq e+1$, we have $n+a \geq n+e+1$ and hence $n-1 \geq n+e-a$. Therefore, $|D_{r_e(n)}| \leq 2^{n-1} < 2^n$ as required. 
\end{proof}
\end{proof}

A natural question about the proof of Theorem \ref{thm:CompTrNotLowIM} is whether this construction could have been performed inside any $\Pi^0_1$ class of 
separating sets. That is, instead of building $A$ and $B$, could we have started with a given pair of computably inseparable sets $U$ and $V$ and defined our path $C$ 
within $\text{Sep}(U,V)$ by a more general argument? To see that the answer is no, consider the class $\text{Sep}(U,V)$ of 
$\{ 0,1 \}$-valued DNC functions. Each $C \in \text{Sep}(U,V)$ computes a \ml random and hence this class of separating sets has no elements which are computably traceable.

\begin{thm}\label{thm:NotCompTrLowIM}
There is a hyperimmune-free degree that is not computably traceable but that is low for isomorphism.
\end{thm}

\begin{proof}
We build a set $A$ with these degree-theoretic properties. To ensure that $A$ is not computably traceable, it suffices to work with the order function $p(n) = 2^n$. 
We need to construct a function $f \leq_T A$ such that for all computable functions $r(n)$ for which $|D_{r(n)}| \leq 2^n$, there is an $n$ such that $f(n) \not \in 
D_{r(n)}$. Let $r_e(n)$, $e \in \omega$, be a (noneffective) listing of all computable functions such that $|D_{r(n)}| \leq 2^n$. 

It will simplify our notation to use a standard coding to move back and forth between a compact subspace of Baire space and Cantor space. Let $T_B$ be the subtree of 
$\omega^{< \omega}$ defined by declaring $\sigma \in T_B$ if and only if for all $n < |\sigma|$, $\sigma(n) \leq 2^n$. The facts that $T_B$ is  
$(2^n+1)$-branching at level $n$ and that $2^n+1 >  |D_{r_e(n)}|$ will enable us to make the degree of $A$ not computably traceable. 

We map $T_B$ onto a subtree of Cantor space by recursively defining an image $\widehat{\sigma} \in 2^{< \omega}$ for each $\sigma \in T_B$. For the 
empty string $\lambda$, let $\widehat{\lambda} = \lambda$. Assume $\widehat{\sigma}$ has been defined and $|\sigma| = n-1$. 
For $i \leq 2^n$, set $\widehat{\sigma*i} = \widehat{\sigma}*0^i1$. Define $T_C$, the image of $T_B$ in Cantor space, to be the downward closure of the 
set of nodes $\widehat{\sigma}$ for $\sigma \in T_B$. We say that $\tau \in T_C$ is a (\textit{level} $n$) \textit{coded node} if $\tau = \widehat{\sigma}$ for some $\sigma \in T_B$ 
(with $|\sigma| = n$). For a coded node $\tau \in T_C$, we let $\sigma_{\tau} \in T_B$ denote the unique node such that $\tau = \widehat{\sigma}_{\tau}$. 

Formally, our construction takes place inside $T_C$, but we often treat coded nodes as though they were elements of $T_B$. For example, if $\tau$ is a coded node, 
then we refer to the nodes $\widehat{\sigma_{\tau}*i}$ as the \textit{coded successors of} $\tau$, and we abuse notation by referring to these nodes as the coded 
successors $\tau*i$ for $i \leq 2^n$. We can effectively determine whether a node is coded or not, and if it is coded, which nodes are its coded successors.

Given a path $A$ in $T_C$, we pull $A$ back to a path $f_A \leq_T A$ in $T_B$ using the initial segments $\tau \subseteq A$, which are coded nodes. 
To define $f_A(n)$, let $\tau \subseteq A$ be a level $m$ coded node for $m > n$ and set $f_A(n) = \sigma_{\tau}(n)$.  

Let $T \subseteq T_C$ be a subtree of $T_C$, and let $\tau$ be a level $n$ coded node such that $\tau \in T$. We say $\tau$ \textit{has a single coded successor in} $T$ if there 
is only one coded successor $\tau*i \in T$. We say $\tau$ \textit{is maximally branching in} $T$ if for every $i \leq 2^n$, the coded successor $\tau*i$ is in $T$. 

We build a sequence of computable subtrees $T_0 = T_C \supseteq T_1 \supseteq \ldots$ such that each tree $T$ in this sequence has the following properties.
\begin{enumerate}
\item[(P1)] Every $\tau \in T$ is extendible to an infinite path in $T$.
\item[(P2)]  For every coded node $\tau \in T$, either $\tau$ has only one coded successor in $T$ or $\tau$ is maximally branching in $T$. The least coded node which is 
maximally branching is called the \textit{pseudoroot of} $T$ and is denoted $\lambda_T$ (not to be confused with the empty string $\lambda$). 
\item[(P3)]  If any level $n$ coded node $\tau \in T$ is maximally branching in $T$, then every level $n$ coded node in $T$ is maximally branching in $T$.
\item[(P4)]  For every $m$, there is an $n > m$ such that every level $n$ coded node in $T$ is maximally branching in $T$.
\end{enumerate}

A computable subtree of $T_C$ satisfying (P1)--(P4) is called an \textit{appropriate subtree}. Since we can effectively determine which coded nodes in an appropriate subtree 
$T$ have a single coded successor and which are maximally branching, it follows from (P4) that we can 
effectively find a sequence of levels $m^T_0 < m^T_1 < \ldots$ such that a coded node $\tau \in T$ is maximally branching if and only if there is an $i$ such that $\tau$ is a 
level $m^T_i$ coded node. By definition, the pseudoroot $\lambda_T$ is a coded node that is an initial segment of every infinite path $A$ in $T$, and furthermore, 
if $\lambda_T$ is a level $m$ coded node and $n < m$, then $f_A(n) = \sigma_{\lambda_T}(n)$. 

With this notation, we can specify our requirements and describe how they can be satisfied using appropriate subtrees. 
To make the degree of $A$ not computably traceable, we meet the requirements 
\[
\mathcal{R}_e: \exists n \left( f_A(n) \not \in D_{r_e(n)} \right).
\]
An appropriate subtree $T$ \textit{satisfies} $\mathcal{R}_e$ if there is $m > n$ such that the pseudoroot $\lambda_T$ is a level $m$ coded node and 
$\sigma_{\lambda_T}(n) \not \in D_{r_e(n)}$. If $A$ is an infinite path in such a tree, then $f_A$ satisfies $\mathcal{R}_e$ because $f_A(n) = \sigma_{\lambda_T}(n) \not \in D_{r_e(n)}$. 

To make $A$ have hyperimmune-free degree, we meet the requirements 
\[
\mathcal{S}_e: \text{If } \Phi_e^A \text{ is total, then there is a computable } g \text{ such that } \forall n \left( \Phi_e^A(n) \leq g(n) \right).
\]
We say that an appropriate subtree $T$ \textit{satisfies} $\mathcal{S}_e$ if one of the following conditions holds.
\begin{enumerate}
\item[(S1)] There is an $n$ such that $\Phi_e^{\tau}(n) \! \uparrow$ for all $\tau \in T$.
\item[(S2)] For every $n$ and every level $m^T_n+1$ coded node $\tau \in T$, $\Phi_e^{\tau}(n) \! \downarrow$.
\end{enumerate}
If $A$ is an infinite path in an appropriate subtree satisfying (S1), then $A$ satisfies $\mathcal{S}_e$ because $\Phi_e^A(n) \! \uparrow$. On the other hand, if $A$ is an infinite 
path in an appropriate subtree $T$ satisfying (S2), then $\Phi_e^A(n) \leq g(n)$ for the computable function $g$ defined by
\[
g(n) = \max \{ \Phi_e^{\tau}(n) \mid \tau \in T \text{ is a level } m^T_n+1 \text{ coded node} \}.
\]

To make the degree of $A$ low for isomorphism, fix a (noneffective) list $(\mathcal{A}_i, \mathcal{B}_i)$, $i \in \omega$, of all pairs of computable directed graphs with 
domain $\omega$. We meet the requirements 
\[
\mathcal{L}_{\langle e,i \rangle}: \text{If } \Phi_e^A \text{ is an isomorphism from } \mathcal{A}_i \text{ to } \mathcal{B}_i, \text{ then } \mathcal{A}_i \cong_{\Delta^0_1} \mathcal{B}_i.
\]
We say that an appropriate subtree $T$ \textit{satisfies} $\mathcal{L}_{\langle e,i \rangle}$ if one of the following conditions holds. 
\begin{enumerate}
\item[(L1)] There is an $n$ such that for all $\tau \in T$, $\Phi_e^{\tau}(n) \! \uparrow$. 
\item[(L2)] There is an $m$ such that for all $\tau \in T$ and $n \in \omega$, if $\Phi_e^{\tau}(n) \! \downarrow$, then $\Phi_e^{\tau}(n) \neq m$.
\item[(L3)] There is an $n$ such that $\Phi_e^{\lambda_T} \!\! \upharpoonright n$ converges but is not a partial graph isomorphism from $\mathcal{A}_i$ to $\mathcal{B}_i$. 
\item[(L4)] For every infinite path $P$ through $T$, $\Phi_e^P$ is a graph isomorphism from $\mathcal{A}_i$ onto $\mathcal{B}_i$. 
\end{enumerate}
Let $A$ be an infinite path through an appropriate subtree $T$ such that $T$ satisfies $\mathcal{L}_{\langle e,i \rangle}$. If $T$ satisfies one of (L1), (L2) or (L3), then $A$ satisfies 
$\mathcal{L}_{\langle e,i \rangle}$ because $\Phi_e^A$ is not total, $\Phi_e^A$ is not onto, or $\Phi_e^A$ is not an isomorphism from $\mathcal{A}_i$ to $\mathcal{B}_i$. 
Suppose $T$ satisfies (L4). Since $T$ is computable and every $\tau \in T$ is extendible, $T$ has an infinite computable path $P$. By (L4), $\Phi_e^P$ is a computable isomorphism 
from $\mathcal{A}_i$ to $\mathcal{B}_i$ and so $\mathcal{L}_{\langle e,i \rangle}$ is satisfied. 

Notice that for each of these types of requirements, if an appropriate subtree $T$ satisfies the requirement and $S \subseteq T$ is an appropriate subtree, then $S$ 
also satisfies the requirement. We are now ready to proceed with the construction. List the requirements in order type $\omega$ with the restriction that for all $e$ and $i$, 
$\mathcal{S}_e$ has higher priority than $\mathcal{L}_{\langle e,i \rangle}$. At stage 0, set $T_0 = T_C$. At stage $s+1$, we define an appropriate subtree $T_{s+1} \subseteq T_s$ 
which satisfies the $s$-th requirement. We now consider the possible forms the $s$-th requirement may take.

For the first case, suppose the $s$-th requirement is $\mathcal{R}_e$. Let $n$ be such that $\lambda_{T_s}$ is a level $n$ coded node. Fix  
$i \leq 2^n$ such that $i \not \in D_{r_e(n)}$ and note that $\lambda_{T_s}$ has a coded successor $\lambda_{T_s}*i \in T_s$ because it is maximally branching in $T_s$. 
Let $T_{s+1} \subseteq T_s$ consist of all nodes in $T_s$ which are comparable with the coded successor $\lambda_{T_s}*i$. $T_{s+1}$ is an appropriate 
subtree satisfying $\mathcal{R}_e$.

For the second case, suppose the $s$-th requirement is $\mathcal{S}_e$. If there is an $n \in \omega$ and a  
$\tau \in T_s$ such that for all nodes $\rho \in T_s$ with $\tau \subseteq \rho$, $\Phi_e^{\rho}(n) \! \uparrow$, then fix such a $\tau$ and let $T_{s+1} \subseteq T_s$ consist of 
all nodes in $T_s$ which are comparable with $\tau$. $T_{s+1}$ is an appropriate subtree satisfying $\mathcal{S}_e$ by (S1). 

If not, then for every $n$ and every coded node $\tau \in T_s$, there is a coded node $\rho \in T_s$ such that $\tau \subseteq \rho$ and 
$\Phi_e^{\rho}(n) \! \downarrow$. We define $T_{s+1}$ as follows. 
Let $k_0$ be such that the pseudoroot $\lambda_{T_s}$ is a level $k_0$ coded node. Add $\lambda_{T_s}$ and its initial segments to $T_{s+1}$.
For each coded successor $\lambda_{T_s}*i$ such that $i \leq 2^{k_0}$, let $\rho_i \in T_s$ be a coded node such that $\lambda_{T_s}*i \subseteq \rho_i$ and 
$\Phi_e^{\rho_i}(0) \! \downarrow$. Without loss of generality, we assume that each $\rho_i$ is a coded node at the same level $k_1$. Add each of the nodes $\rho_i$ 
and their initial segments to $T_{s+1}$. 

Now, repeat this process to construct the rest of $T_{s+1}$ recursively. For each level $k_1$ coded node $\rho_i$, consider each coded successor 
$\rho_i*j$, $j \leq 2^{k_1}$, and fix a coded node $\rho_{i,j} \in T_s$ such that $\rho_i*j \subseteq \rho_{i,j}$ and $\Phi_e^{\rho_{i,j}}(1) \! \downarrow$. Without loss of 
generality, each $\rho_{i,j}$ is a coded node at the same level $k_2$. Add each of the nodes $\rho_{i,j}$ and their initial segments to $T_{s+1}$. We continue in this fashion to 
define levels $k_0 < k_1 < \ldots$ and level $k_{\ell +1}$ coding nodes $\rho_{i_0, \ldots, i_\ell}$ such that $\Phi_e^{\rho_{i_0, \ldots, i_\ell}}(n) \! \downarrow$ for $n \leq \ell$. 
By construction, a coded node $\tau \in T_{s+1}$ is maximally branching in $T_{s+1}$ if $\tau$ is a level $k_{\ell}$ node for some $\ell$ and otherwise $\tau$ has a single 
coded successor. It follows that $T_{s+1}$ is an appropriate subtree satisfying $\mathcal{S}_e$ by (S2). 

For the last case, suppose the $s$-th requirement is $\mathcal{L}_{\langle e,i \rangle}$. By the ordering of requirements, the appropriate subtree $T_s$ satisfies $\mathcal{S}_e$. We break into 
subcases. First, suppose $T_s$ satisfies (S1). Let $T_{s+1} = T_s$ and note that $T_{s+1}$ satisfies (L1). For the remaining subcases, we assume $T_s$ satisfies (S2). 
Second, suppose there is an $n$ and a level $m_n^{T_s}$ coded node $\tau \in T_s$ such that $\Phi_e^{\tau} \upharpoonright n$ is not a partial graph isomorphism from 
$\mathcal{A}_i$ to $\mathcal{B}_i$. Fix such a node $\tau$ and let $T_{s+1}$ be the set of nodes in $T_s$ which are comparable with $\tau$. In this case, $T_{s+1}$ is an 
appropriate subtree satisfying (L3). Third, suppose there is an $m$ and a coded node $\tau \in T_s$ such that for all $\rho \in T_s$ with $\tau \subseteq \rho$, 
if $\Phi_e^{\rho}(n) \! \downarrow$, then $\Phi_e^{\rho}(n) \neq m$. Fix such a node $\tau$ and let $T_{s+1}$ be the set of nodes in 
$T_s$ which are comparable with $\tau$. In this case, $T_{s+1}$ is an appropriate subtree satisfying (L2). 

Finally, if this case also fails, then we know that for every coded node $\tau \in T_s$ and every $m$, there is a coded node $\rho \in T_s$ and an $n$ such that $\tau \subseteq \rho$ 
and $\Phi_e^{\rho}(n) = m$. By a process similar to the one used to define an appropriate tree satisfying (S2) above, there is an appropriate subtree $T_{s+1} \subseteq T_s$ 
such that for every coding node $\rho$ extending the $k$-th maximally branching level $m_k^{T_{s+1}}$ in $T_{s+1}$, $k \in \text{range}(\Phi_e^{\rho})$. Therefore, for every 
infinite path $P$ in $T_{s+1}$ and every $m$, $\Phi_e^P$ is total and is onto $\omega$. Since for each $n$ and each level $m_n^{T_{s+1}}$ coded node $\tau \in T_{s+1}$, 
$\Phi_e^{\tau} \upharpoonright n$ is a partial isomorphism from $\mathcal{A}_i$ to $\mathcal{B}_i$, it follows that $T_{s+1}$ is an appropriate subtree satisfying (L4). 

To define the set $A$, note that by construction, $\lambda_{T_s} \subseteq \lambda_{T_{s+1}}$ and the lengths of the pseudoroots are unbounded. Let $A$ be the unique set 
such that for every $s$, $\lambda_{T_s}$ is an initial segment of $A$. Since $A$ is an infinite path on each tree $T_s$, $A$ satisfies all of the requirements.
\end{proof}

\def\cprime{$'$}

\end{document}